\documentclass[a4paper,12pt,reqno]{amsart}

\usepackage{amsmath,amssymb,amsthm,latexsym}

\newtheorem{theorem}{Theorem}
\newtheorem{corollary}{Corollary}
\newtheorem{definition}{Definition}

\begin{document}

\title[Another Monotone Selection Principle]%
{A Non-commutative Monotone Selection Principle}

\author[M. Thill]{Marco Thill}

\address{bd G.\,-\,D.\ Charlotte 53 \\
                 L\ -\,1331 Luxembourg\,-\,City}

\email{math@pt.lu}

\date{}



\thanks{My thanks go to Michael Frank for the impetus necessary for this paper.}

\begin{abstract}
We give an elementary proof of a monotone selection principle which allows
to pass from increasing nets to increasing sequences in the Hermitian part
of a $\sigma$-finite von Neumann algebra. This is to be seen as a
``monotone version'' of first countability.
\end{abstract}

\maketitle


This compuscript is a sequel to the previous arXived article of the author \cite{Thi},
where similar results were obtained in the commutative case.

\begin{definition}
Let $W$ be an ordered vector space, and let $V$ be a vector subspace of $W$.
(Typically $W$ the Hermitian part of a C*-algebra.)

We shall say that $W$ is \emph{monotone complete}, if each non-empty upper
bounded upward directed subset of $W$ has a supremum in $W$. We shall say that $W$
is \emph{monotone sequentially complete}, if each upper bounded increasing sequence
in $W$ has a supremum in $W$.

If $W$ is monotone complete, then $V$ shall be called \emph{monotone closed}
in $W$, if $V$ contains the supremum $\sup J$ in $W$ of each non-empty upward directed
subset $J$ of $V$ that is upper bounded in $W$. (It is clear that then $\sup J$ also is
the supremum in $V$ of $J$, and that $V$ then is monotone complete.) If $W$ is monotone
sequentially complete, then $V$ shall be called \emph{monotone sequentially closed} in
$W$, if $V$ contains the supremum $\sup _n j_n$ in $W$ of each increasing sequence
$(j_n)$ in $V$ that is upper bounded in $W$. (It is clear that then $\sup _n j_n$ also is the
supremum in $V$ of $(j_n)$, and that $V$ then is monotone sequentially complete.)

If $W _+$ denotes the set of positive elements of $W$, then a positive linear functional
$\psi$ on $W$ shall be called \emph{faithful}, if $\psi (a) > 0$ for each
$a \in W _+ \setminus \{ 0 \}.$
\end{definition}

For example the Hermitian part a von Neumann algebra is monotone complete.
A C*-algebra of bounded operators on a Hilbert space $H$ is a von Neumann algebra
if and only if its Hermitian part is monotone closed in the Hermitian part of $B(H)$.
For a proof of this theorem of Kadison, see for example \cite[2.2.4]{Ped}.

A C*-algebra of bounded operators on a Hilbert space $H$ is called a Borel
$\ast$-algebra if its Hermitian part is monotone sequentially closed in the Hermitian
part of $B(H)$. (See \cite[4.5.5]{Ped}.)

\begin{definition}
Let $W_1$, $W_2$ be ordered vector spaces and let $\phi$ be a positive
linear map from $W_1$ to $W_2$.

If $W_1$, $W_2$ are monotone complete, then $\phi$ is called \emph{normal}, if
for each non-empty upper bounded upward directed subset $J$ of $W_1$, one has
$\phi (\sup J) = \sup _{j \in J} \phi (j)$. If $W_1$, $W_2$ are monotone sequentially
complete, then $\phi$ is called \emph{sequentially normal}, if for each upper bounded
increasing sequence $(j_n)$ in $W_1$, one has $\phi (\sup _n j_n) = \sup _n \phi (j_n)$.
\end{definition}

A von Neumann algebra is $\sigma$-finite if and only if its Hermitian part carries
a faithful normal positive linear functional. (See e.g.\ \cite[2.5.6]{BR}.)

\begin{theorem}
Let $W$ be an ordered vector space. Assume that $W$ is monotone complete
and carries a faithful normal positive linear functional $\psi$.
(E.g.\ $W$ the Hermitian part of a $\sigma$-finite von Neumann algebra.)

Then the following \emph{monotone selection principle} holds. Whenever $J$
is a non-empty upper bounded upward directed subset of $W$, there exists
an increasing sequence $(j_n)$ in $J$ such that $\sup _n j_n = \sup J$. Every
increasing sequence $(j_n)$ in $J$ with
\[ \sup _{n \geq 1} \psi (j_n) = \sup _{j \in J} \psi (j) \]
does the job.
\end{theorem}

\begin{proof}
Let $J$ be a non-empty upper bounded upward directed subset of $W$.
Let $(j_n) _{n \geq 1}$ be any increasing sequence in $J$ with
\[ \sup _{n \geq 1} \psi (j_n) = \sup _{j \in J} \psi (j). \tag*{$(*)$} \]
With $s := \sup J \in W$, we have
\[ \psi(s) = \sup _{j \in J} \psi (j) \tag*{$(**)$} \]
by normality of $\psi$.
We define $j_0 := \sup _{n \geq 1} j_n \in W$. Then $\psi (j_0) = \psi (s)$ by
sequential normality of $\psi$ and the equations $(*)$ and $(**)$. On one hand,
we then have $s - j_0 \geq 0$. On the other hand, we have $\psi (s - j_0) = 0$.
Therefore, the fact that $\psi$ is faithful implies that $s = j_0$.
\end{proof}

\begin{corollary}
Let $W$ be an ordered vector space. Assume that $W$ is monotone complete
and carries a faithful normal positive linear functional $\psi$. Then every monotone
sequentially closed subspace of $W$ is monotone closed in $W$.
\end{corollary}

We thus have an elementary proof of the following well-known result,
cf.\ \cite[4.5.5]{Ped}.

\begin{corollary}
A Borel $\ast$-subalgebra of a $\sigma$-finite von Neumann algebra
is itself a von Neumann algebra.
\end{corollary}

\begin{theorem}
Let $W_1$, $W_2$ be monotone complete ordered vector spaces. Assume
that $W_1$ carries a faithful normal positive linear functional. Then every
sequentially normal positive linear map from $W_1$ to $W_2$ is normal.
\end{theorem}

\begin{corollary}
Let $A$, $B$ be von Neumann algebras, and let $\phi : A \to B$ be a
positive linear map. Assume that $A$ is $\sigma$-finite. If $\phi$ is
sequentially normal, then it is normal.
\end{corollary}

\end{document}